\documentclass{amsart}
\usepackage{graphicx}
\usepackage{enumerate}
\usepackage{fullpage}
\usepackage{subfigure}
\newtheorem{theorem}{Theorem}[section]

\theoremstyle{definition}

\theoremstyle{remark}

\numberwithin{equation}{section}

\newcommand{\set}[1]{\left\{#1\right\}}
\newcommand{\comb}[2]{\left( \begin{array}{c} #1 \\ #2 \end{array} \right)}
\newcommand{\R}{\mathbb R}
\newcommand{\C}{\mathbb C}

\newcommand{\BB}{{\mathcal B}}
\newcommand{\LL}{{\mathcal L}}
\newcommand{\PP}{{\mathcal P}}
\newcommand{\HH}{{\mathcal H}}

\newcommand{\eps}{\varepsilon}
\newcommand{\hyper}[2]{{}_{_{#1}}F_{_{#2}}}
\begin{document}

\title[]{Escape rates for the Farey map with approximated holes}%
\author{C. Bonanno} 
\address{Dipartimento di Matematica, Universit\`a di Pisa, Largo Bruno Pontecorvo 5, 56127 Pisa, Italy}
\email{claudio.bonanno@unipi.it}
\author{I. Chouari}
\address{Faculty of Science, Computational Mathematics Laboratory, University of Monastir, Monastir 5000, Tunisia}
\thanks{CB would like to thanks the organizers of the workshop ``DinAmicI IV'', the meeting of the group of Italian researchers in dynamical systems (see http://www.dinamici.org/), because at the workshop he heard about \cite{KM}, which has inspired this paper. IC acknowledges warm hospitality and financial support from the Dipartimento di Matematica, Universit\`a di Pisa, where part of this paper has been written.}
\maketitle

\begin{abstract}
We study the escape rate for the Farey map, an infinite measure preserving system, with a hole including the indifferent fixed point. Due to the ergodic properties of the map, the standard theoretical approaches to this problem cannot be applied. To overcome this difficulties we propose here to consider approximations of the hole by means of real analytic functions. We introduce a particular family of approximations and study numerically the behavior of the escape rate for ``shrinking'' approximated holes. The results suggest that the scaling of the escape rate depends on the chosen approximation, but ``converges'' to the behavior found for piecewise linear approximations of the map in \cite{KM}.
\end{abstract}


\section{Introduction} 

In recent years there has been a quick growth of the number of papers dealing with statistical properties of dynamical systems with holes. The origin of these studies can be found in the paper \cite{PY}, where it was posed the question about the statistical properties of the motion of a particle inside a billiard table with a small hole.  For example, if $p_n$ is the probability that a trajectory remains on the table until time $n$, what is the decay rate of $p_n$? In general, they asked whether an initial distribution would converge, under suitable renormalization, to some limit distribution, a conditionally invariant measure. Much attention to these problems has been paid also by the physics community, see e.g. \cite{Alt} and references therein.

In a general dynamical system, the first question posed in \cite{PY} has become the following: let $X$ be the phase space of a dynamical system with $m$ a distribution of initial conditions; let $H$ be a hole in the space $X$; let $p_n:= m(S_n)$, where $S_n$ is the set of surviving points up to time $n$; study the decay rate
\[
\gamma := \lim_{n\to \infty}\, - \frac{\log(p_n)}{n}\, ,
\]
which is called the \emph{escape rate}. The existence of the limit and its dependence on the initial distribution $m$ has been studied in \cite{DY}, where it is shown that in the ideal case it is possible to compute the escape rate by using the transfer operator associated to the system. 

Let $X=[0,1]$ and $F:X\to X$ be a smooth map with a finite number of pre-images for each $x\in X$, then the \emph{transfer operator} associated to $F$ is defined as
\[
(\PP f)(x) := \sum_{y\in F^{-1}(x)}\, \frac{f(y)}{|F'(y)|}\, .
\]
The operator $\PP$ has spectral radius less or equal than 1, and if there exists a function $g$ such that $\PP g = g$, then $d\mu(x) = g(x) dx$ is an $F$-invariant measure. We refer to \cite{Ba} for more properties of the transfer operators. When there is a hole $H$ in $X$, one can consider the \emph{transfer operator for the open system} 
\[
\PP^{op} f := \PP ((1-\chi_{H})\, f)\, ,
\]
where $\chi_H$ is the indicatrix function of the set $H$. Namely the transfer operator for the open system considers only pre-images of a point $x\in X$ which are not in the hole $H$. Then the escape rate is obtained by
\begin{equation}\label{escape}
\gamma = -\log(\lambda_H)
\end{equation}
where $\lambda_H$ is the largest eigenvalue of $\PP^{op}$.

The escape rate has been largely studied mainly for hyperbolic systems and piecewise expanding maps of the interval. We refer to \cite{DY} for references. In this paper we are interested in studying the behavior of the escape rate as the hole shrinks. Rigorous results on this asymptotic behavior are given in \cite{LK} for piecewise expanding maps of the interval, for which it is found that if $|H|=\epsilon$ then $\gamma \sim const. \times \epsilon$ as $\epsilon \to 0^+$. Higher order corrections can be found in \cite{Cri,D}.

However, it is well known that statistical properties of dynamical systems dramatically change passing from uniformly hyperbolic systems to intermittent ones, that is systems which have an indifferent fixed point. And in particular when the intermittent system preserves only one absolutely continuous invariant measure which is infinite. Intermittent systems have been introduced in the mathematical physics literature in \cite{PM} as a simple model of the phenomenon of intermittency, that is the alternation of a turbulent and a laminar phase in a fluid. As dynamical systems on the unit interval $[0,1]$, they may be represented by the family of maps $F(x) = x+ x^\alpha$ (mod 1), with $\alpha>1$, which have a fixed point at $x=0$ with $F'(0) = 1$ and $F'(x)- 1 \approx x^{\alpha-1}$ as $x\to 0^+$. For $\alpha\in (1,2)$ the system has a finite absolutely continuous invariant measure. For this case, the escape rate has been recently studied in \cite{DF}, where it is shown that the probability $p_n$ may decrease polynomially, so in the definition of $\gamma$ one should divide by $log (n)$. However the results in \cite{DF} consider the case of a hole $H$ which is generated by the Markov partition of the map, and such that the indifferent fixed point is outside $H$. Hence the polynomial escape rate is a consequence of the typical slower decay of correlations found for intermittent systems.

\begin{figure}[h]
\begin{center}
\includegraphics[width=8cm]{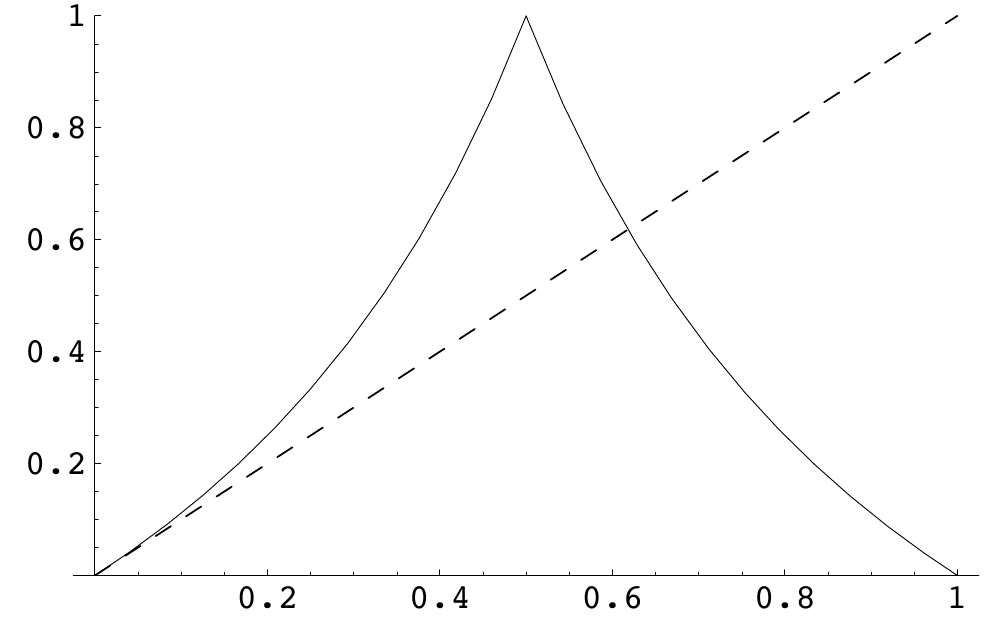}
\caption{The Farey map.} \label{gr-farey}
\end{center}
\end{figure}

In this paper we are interested in the case of an intermittent dynamical system with infinite absolutely continuous invariant measure. This corresponds to values $\alpha \ge 2$ in the example $F(x) = x+ x^\alpha$ (mod 1). This case has been recently studied by G. Knight and S. Munday in \cite{KM}. They study the asymptotic behavior of the escape rate for vanishing hole and prove that different behaviors are possible. Their methods are analytical and use the definition of the \emph{dynamical zeta function} associated to a system, and in particular the relations with the transfer operator. However, to apply these methods they have to consider piecewise linear maps and holes which are generated by the partition associated to the map. To give a flavor of their results, we consider a typical example of intermittent dynamical system on the interval $[0,1]$ with infinite measure, the \emph{Farey map}. The Farey map is defined by
\begin{equation} \label{farey}
F(x)=\left\{
\begin{array}{ll}
\frac{x}{1-x}\, , & \mbox{if }\ 0\le x\le \frac{1}{2}\\[0.3cm]
\frac{1-x}{x}\, , & \mbox{if }\ \frac{1}{2} \le x \le 1
\end{array} \right.
\end{equation}
and is studied in particular for its relations with the continued fractions expansion of real numbers (see e.g. \cite{BI}). Its graph is shown in Figure \ref{gr-farey}.
In \cite{KM}, the authors consider a piecewise version of $F$, which is obtained by considering the partition $A=\{ (\frac{1}{n+1}, \frac 1n)\}$ and defining
\[
F_p(x)=\left\{
\begin{array}{ll}
2-2x\, , & \mbox{if }\ \frac{1}{2} \le x \le 1\\[0.3cm]
\frac{n+1}{n-1}\, x - \frac{1}{n(n-1)}\, , & \mbox{if }\ \frac{1}{n+1} \le x \le \frac{1}{n}
\end{array} \right.
\]
Then for holes $H_n = [0,\frac 1n)$, they show that
\begin{equation}\label{kn-mu}
\gamma \approx \frac{1}{n\, \log n} = \frac{|H_n|}{-\log |H_n|}\qquad \text{as}\quad n\to \infty\, ,
\end{equation}
which is a slightly faster decay than those proved for expanding maps in \cite{LK}. More general behavior are found in \cite{KM} by varying the partition $A=\{(a_{n+1},a_n)\}$ with $a_n\to 0^+$, but always choosing a hole of the form $[0,a_n)$.

The aim of this paper is to study the escape rate for the Farey map \eqref{farey} with hole $H=[0,\epsilon)$, and its asymptotic behavior as $|H| = \eps\to 0^+$. To our knowledge this is the first case where this problem is studied in such a generality for an infinite measure preserving dynamical system. We use \eqref{escape} to compute the escape rate through the transfer operator $\PP^{op}$ associated to $F$ with a hole $H$. The main problem is that due to the presence of the indifferent fixed point, the transfer operator of the Farey map doesn't have the so-called \emph{spectral gap}, that is a gap between the maximal eigenvalue and the essential spectrum, independently of the smoothness of functions on which the operator is applied (see \cite{CI}). The spectral gap is found for example in the study of the transfer operator of expanding maps, when restricted to the space of bounded variation functions, and is fundamental to obtain exponential decay of correlations for the map and is used also in \cite{LK} to study the asymptotic behavior of the escape rate.

A good space of functions $\HH$ on which to study the spectral properties of the transfer operator of the Farey map has been introduced in \cite{Is,BGI}. It consists of holomorphic functions on the disc $\{|z-\frac 12|< \frac 12\}$ obtained as integral transform of functions on the positive real axis (see \eqref{spazio} below). So we should consider the action of $\PP^{op}$ on $\HH$. Unfortunately, since the indicatrix function $\chi_H$ is not real analytic, we have to introduce an approximation of $\chi_H$. Hence we study the escape rate for the Farey map \eqref{farey} with approximated holes which contain the indifferent fixed point, and are general in the sense that we don't require any relation between the length of the hole and the map.

In Section \ref{sec:to} we introduce the operator $\PP^{op}$ and its approximation $\tilde \PP^{op}$, obtained by approximating the function $\chi_H$ with a real analytic function. In Section \ref{sec:ma} we introduce the approach we use to study the spectral properties of $\tilde \PP^{op}$, and show our results. The conclusions of our study are exposed in the final section.

\section{Transfer operators} \label{sec:to}

\subsection{The Farey map}

The transfer operator $\mathcal{P}$ associated to the map $F$ acts on functions $f:[0,1]\to \C$ as
\[
(\mathcal{P} f)(x) := \sum_{y\, :\, F(y)=x}\, \frac{f(y)}{|F'(y)|}
\]
which using \eqref{farey} becomes
\[
(\mathcal{P} f)(x)=(\mathcal{P}_{0}f+\mathcal{P}_{1}f)(x)
\]
with
\begin{equation}  \label{p0p1}
(\mathcal{P}_{0}f)(x)=\frac{1}{(1+x)^{2}}\, f\Big(\frac{x}{1+x}\Big)\quad \text{and} \quad (\mathcal{P}_{1}f)(x)=\frac{1}{(1+x)^{2}}\, f\Big(\frac{1}{1+x}\Big)\, .
\end{equation} 

The operator $\mathcal{P}$ has been studied in \cite{BGI,pap1} on the invariant Hilbert space $\HH$ defined as
\begin{equation}\label{spazio}
\HH:= \set{f:[0,1]\to \C\, :\, f=\mathcal{B}[\varphi]\ \text{for some } \varphi\in L^{2}(m)}
\end{equation}
where $\mathcal{B}[\cdot]$ denotes the generalized Borel transform
\begin{equation} \label{trans}
(\mathcal{B}[\varphi])(x):= \frac{1}{x^{2}}\int_{0}^{\infty}e^{-\frac{t}{x}}\, e^{t}\, \varphi(t)\, dm(t)\, ,
\end{equation}
and $L^{2}(m):= L^2(\R^+,m)$ where $m$ is the measure on $\R^+$
\[
dm(t)=te^{-t}dt.
\]
The space $\HH$ is endowed with the inner product inherited by the inner product on $L^2(m)$ through the $\mathcal{B}$-transform, that is
\[
(f_{1},f_{2})_{\HH}:=\int_{0}^{\infty}\varphi_{1}(t)\, \overline{\varphi_{2}(t)}\,dm(t)\qquad \text{if}\quad f_{i}=\mathcal{B}[\varphi_{i}]\, .
\]
It is immediate to show that functions in $\HH$ can be continued to holomorphic functions on the disc $\{|z-\frac 12|<\frac 12\}$, and moreover eigenfunctions of $\PP$ in $\HH$ can be continued to holomorphic function on the half-plane $\{\Re(z)>0\}$. From the computational point of view, it is more convenient to use the $\mathcal{B}$-transform to read the action of $\PP$ on $L^2(m)$. Indeed
\begin{equation}\label{traduzione}
\mathcal{P} \, \mathcal{B}[\varphi]=\mathcal{B}[(M+N)\varphi]
\end{equation}
for all $\varphi \in L^2(m)$, where $M,N :L^2(m) \rightarrow L^2(m)$ are self-adjoint bounded linear operators defined by
\begin{equation}\label{trad}
(M \varphi)(t)= e^{-t}\, \varphi(t) \qquad \text{and} \qquad (N\varphi)(t)= \int_{0}^{\infty}J_{1}\left(2\sqrt{st}\right)\, \sqrt{\frac{1}{st}}\ \varphi(s)\, dm(s)
\end{equation}
where $J_{q}$ denotes the Bessel function of order $q$. 

The spectrum of the operator $\PP$ is the unit interval $[0,1]$, and it is continuous spectrum without eigenvalues on $\HH$. This phenomenon is due to the presence of the indifferent fixed point at the origin. In particular $1$ is not an eigenvalue, and indeed the Farey map has an infinite absolutely continuous invariant measure with density $\frac 1x$, which does not belong to $\HH$.

\subsection{The Farey map with a hole}

\begin{figure}[ht]
\begin{center}
\includegraphics[width=8cm,keepaspectratio]{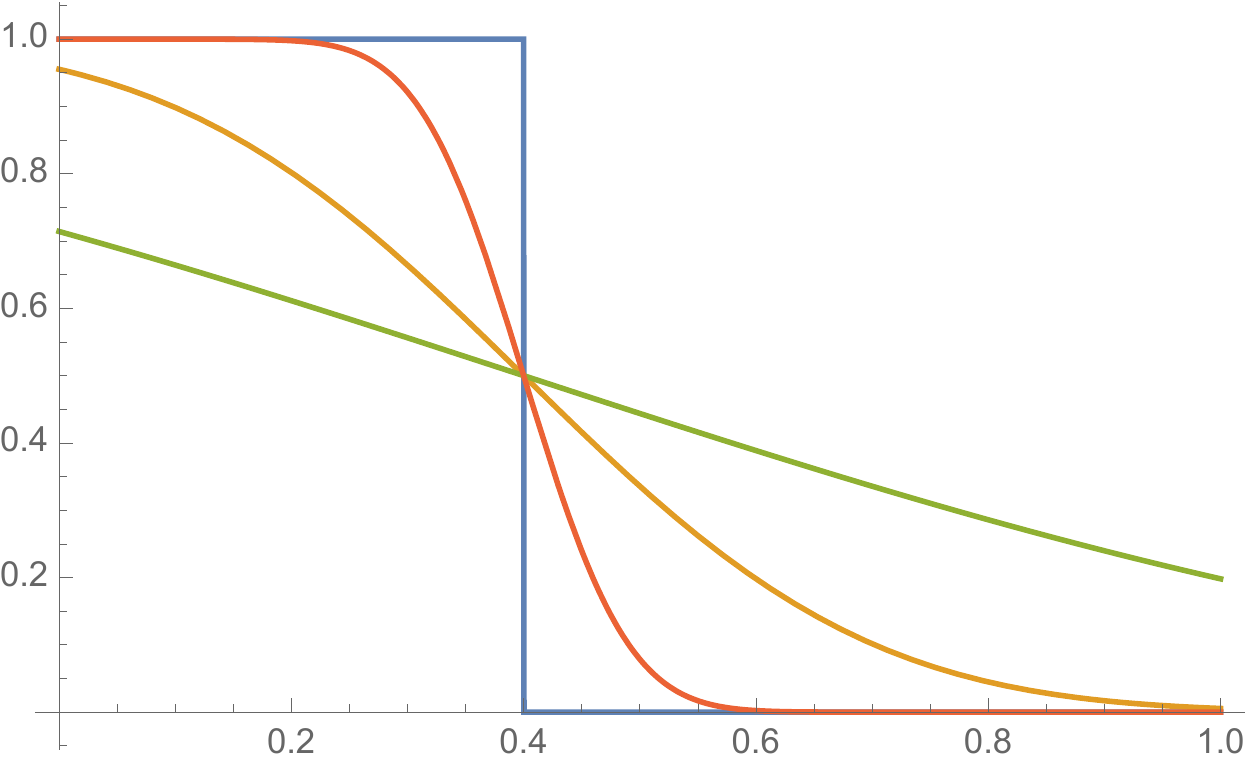}
\end{center}
\caption{The graph of the function $\xi_\mu(x,a)$ for $a=0.4$ and $\mu = 1,3,10$, and of the characteristic function of $[0,a]$.} \label{graph}
\end{figure}

The transfer operator $\PP^{op}$ for the map $F$ with a hole in $H=[0,\epsilon)$ is obtained by $\PP$ simply subtracting the contribution from the hole, that is
\[
\PP^{op} f := \PP ((1-\chi_H) f) = \PP f - \PP (\chi_H f)\, ,
\]
where $\chi_H(x)$ is the indicator function of the set $H$. First of all we notice that for $\epsilon<\frac 12$ it holds
\[
\PP(\chi_H f) = \PP_{0} (\chi_H f)\, ,
\]
since $\frac{1}{1+x}\ge \frac 12$ for all $x\in [0,1]$, hence by \eqref{p0p1}
\[
(\PP (\chi_H f))(x) = \frac{1}{(1+x)^{2}}\, \chi_H\Big(\frac{x}{1+x} \Big)\, f\Big(\frac{x}{1+x}\Big) = \frac{1}{(1+x)^{2}}\, \chi_{\tilde H}(x)\, f\Big(\frac{x}{1+x}\Big) = \chi_{\tilde H}(x)\, (\PP_{0} f)(x)
\]
where $\tilde H = [0,\frac{\epsilon}{1-\epsilon})$. Hence
\begin{equation}\label{to-hole-d}
\PP^{op} f = \PP f - \chi_{\tilde H}\, \PP_{0} f\, .
\end{equation}
However, to study the spectral properties of the operator $\PP^{op}$ on $\HH$, we need an approximation of the indicator function since we are dealing with real analytic functions. We consider the family of approximations $\set{ \xi_\mu(x,a)}_{\mu \in\R}$ for the indicator function of an interval $[0,a)$ given by
\begin{equation} \label{famiglia}
\xi_\mu(x,a):=  \frac 12 - \frac 12\text{Erf}\left(\mu \, \left(x  - a\right)\right)\, ,
\end{equation}
where $\text{Erf}$ is the ``error function"
\begin{equation} \label{series}
\text{Erf}(x) := \frac{2}{\sqrt{\pi}}\, \int_0^x\, e^{-t^2}\, dt = \frac{2}{\sqrt{\pi}}\, \sum_{n=0}^\infty\, \frac{(-1)^n \, x^{2n+1}}{n!\, (2n+1)}\, .
\end{equation}
In Figure \ref{graph} we have plotted the function $\xi_\mu(x,a)$ for $a=0.4$ and $\mu=1,3$ and $10$, against $\chi_{[0,a)}(x)$. The functions $\xi_\mu(x,a)$ converge to $\chi_{[0,a)}(x)$ as $\mu\to \infty$ for all $x\not= a$, but clearly not uniformly.

Hence we define a family of approximated transfer operators for the Farey map with hole $[0,\epsilon)$ by using $\xi_\mu$ in \eqref{to-hole-d}, and let
\begin{equation}\label{to-hole-new}
(\tilde \PP^{op}_\mu f)(x) := (\PP f)(x) - \xi_\mu\left(x, \frac{\epsilon}{1-\epsilon}\right)\, (\PP_{0} f)(x)\, .
\end{equation}
Moreover, to study the action of $\tilde \PP^{op}_\mu$ on $\HH$ as for the action of $\PP$, we need the equivalent of equation \eqref{traduzione} for $\tilde \PP^{op}_\mu$. To this aim, we first recall that by definition \eqref{trans}
\[
\BB[\varphi(t)](x) = \frac{1}{x^2}\, \LL[t\, \varphi(t)]\left( \frac 1x\right)
\]
where $\LL[\cdot]$ denotes the standard Laplace transform, then for all $\psi,\varphi \in L^2(m)$ we can write
\begin{equation} \label{serve}
\begin{aligned}
& \BB[\psi](x) = \xi_\mu(x,a) \, \BB[\varphi](x) \quad \Leftrightarrow \quad \LL[t\, \psi(t)]\left(\frac 1x\right) = \xi_\mu(x,a) \,  \LL[t\, \varphi(t)]\left(\frac 1x\right) \quad \Leftrightarrow\\
& \Leftrightarrow \quad \LL[t\, \psi(t)](x) = \xi_\mu\left(\frac 1x, a\right)\, \LL[t\, \varphi(t)](x) \quad \Leftrightarrow \quad t\, \psi(t) = \int_0^t\, s\, \varphi(s)\, \LL^{-1}\left[\xi_\mu\left(\frac 1x, a\right)\right](t-s)\, ds
\end{aligned}
\end{equation}
where in the last equivalence we have used the standard properties for the Laplace transform of convolution of functions. Using \eqref{traduzione} and \eqref{serve} in \eqref{to-hole-new}, we obtain
\begin{equation}\label{traduzione-open-gen}
\begin{aligned}
& (\tilde \PP^{op}_\mu\, \BB[\varphi])(x) = \BB[(M+N)\varphi](x) - \xi_\mu\left(x, \frac{\epsilon}{1-\epsilon}\right)\, \BB[M \varphi](x) = \\
& = \BB[(M+N)\varphi](x) - \BB\left[\frac 1t\, \int_0^t\, s\, (M\varphi)(s)\, \LL^{-1}\left[\xi_\mu\left(\frac 1x, \frac{\epsilon}{1-\epsilon}\right)\right](t-s)\, ds\right] (x) = \\
& = \BB[(\tilde M_\mu + N) \varphi](x)
\end{aligned}
\end{equation}
with $M,N$ as in \eqref{trad}, and 
\begin{equation} \label{emme-tilde}
(\tilde M_\mu \varphi) (t) := (M\varphi)(t) - \frac 1t\, \int_0^t\, s\, (M\varphi)(s)\, \LL^{-1}\left[\xi_\mu\left(\frac 1x, \frac{\epsilon}{1-\epsilon}\right)\right](t-s)\, ds\, .
\end{equation}

\begin{theorem}\label{funzionale}
For $\xi_\mu(x,a)$ as in \eqref{famiglia}, the operator $\tilde M_\mu$ is given by
\begin{equation} \label{trad-op}
\begin{aligned}
(\tilde M_\mu \varphi) (t) & = \frac 12 \Big(1 - \mathrm{Erf}\Big(\frac{\mu\, \epsilon}{1-\epsilon}\Big)\Big)\, (M\varphi)(t) + \\[0.2cm]
& + \frac{1}{\sqrt{\pi}}\, \sum_{n=0}^\infty\, \sum_{k=1}^{2n+1}\, \comb{2n+1}{k}\, \frac{(-1)^{n+k-1}\, \mu^{2n+1}}{n!\, (k-1)!\, (2n+1)}\, \left( \frac{\epsilon}{1-\epsilon} \right)^{2n+1-k}\ \frac 1t\, \int_0^t\, s (M\varphi)(s)\, (t-s)^{k-1}\, ds
\end{aligned}
\end{equation}
and it is a bounded operator $\tilde M_\mu :L^2(m) \to L^2(m)$.
\end{theorem}

\section{The matrix approach and the numerical results} \label{sec:ma}

As shown in \cite{BGI}, the Hilbert space $L^2(m)$ admits a complete orthogonal system $\{e_\nu\}_{\nu\ge 0}$ given by the Laguerre polynomials defined as
\begin{equation} \label{laguerre}
e_\nu(t) := \sum_{m=0}^\nu\, \comb{\nu+1}{\nu-m}\, \frac{(-t)^m}{m!}
\end{equation}
which satisfy 
\[
(e_{\nu},e_{\nu})=\frac{\Gamma(\nu+2)}{\nu!} = \nu+1
\]
for all $\nu\ge 0$. Hence, using \eqref{traduzione-open-gen}, we can study the action of $\tilde \PP^{op}_\mu$ on $\HH$ by the action on $L^2(m)$ of an infinite matrix representing the operators $\tilde P^{op}_\mu:=\tilde M_\mu+N$, defined in \eqref{trad} and \eqref{trad-op}, for the basis $\{e_\nu\}_{\nu\ge 0}$. That is for any $\phi\in L^{2}(m)$, we can write
\[
\phi(t)=\sum_{\nu=0}^{\infty}\phi_{\nu}e_{\nu}(t) \quad \text{with} \quad \phi_\nu = \frac{1}{\nu+1}\, (\phi,e_{\nu})
\] 
hence $\phi$ is an eigenfunction of $\tilde P^{op}_\mu$ with eigenvalue $\lambda$ if and only if
\[
(\tilde P^{op}_\mu \phi,e_{j})=\lambda\, (\phi,e_{j})=\lambda \, (j+1)\, \phi_{j} \qquad \forall \; j\geq 0
\]
Using the notation $c_{j\nu}^\mu := (\tilde P^{op}_\mu e_{\nu},e_{j})$ we obtain that 
\begin{equation} \label{eig-mat}
\tilde P^{op}_\mu \phi=\lambda \phi \quad \Leftrightarrow\quad C_{\mu}\phi=\lambda D\phi \quad \Leftrightarrow\quad A_{\mu}\phi=\lambda \phi
\end{equation}
where $C_{\mu}$ and $D$ are given by 
\[
C_{r}=(c_{j\nu}^{\mu})_{j,\nu\geq 0}\quad \text{and} \quad D=\text{diag}(j+1)_{j\ge 0}
\]
and $A_{\mu}$ is the infinite matrix
\begin{equation}\label{ar}
A_{\mu}= (a_{j\nu}^{\mu})_{j,\nu\geq 0}\qquad \text{with} \quad a_{j\nu}^\mu = \frac{c_{j\nu}^{\mu}}{j+1}\, .
\end{equation}
We have then to compute the terms $c_{j\nu}^{\mu}:= (\tilde P^{op}_r e_{\nu},e_{j})$. From \cite[Prop. 3.1]{pap2} we have
\[
\begin{aligned}
&\frac{1}{j+1}\, (M e_{\nu},e_{j}) = \comb{\nu+j+1}{\nu}\, \frac{1}{2^{\nu+j+2}} \\[0.2cm]
&\frac{1}{j+1}\, (N e_{\nu},e_{j}) = \sum_{\ell=0}^{\nu}(-1)^{\ell}\, \comb{\nu+1}{\nu-\ell}\, \comb{\ell+j+1}{l}\, \frac{1}{2^{\ell+j+2}}\, .
\end{aligned}
\]
Hence
\begin{equation}\label{final}
\begin{aligned}
& a^\mu_{j\nu} = \Big(1 - \mathrm{Erf}\Big(\frac{\mu\, \epsilon}{1-\epsilon}\Big)\Big)\, \comb{\nu+j+1}{\nu}\, \frac{1}{2^{\nu+j+3}} + \sum_{\ell=0}^{\nu}(-1)^{\ell}\, \comb{\nu+1}{\nu-\ell}\, \comb{\ell+j+1}{l}\, \frac{1}{2^{\ell+j+2}} +\\[0.2cm]
& + \frac{1}{(j+1)\sqrt{\pi}}\, \sum_{n=0}^\infty\, \sum_{k=1}^{2n+1}\, \comb{2n+1}{k}\, \frac{(-1)^{n+k-1}\, \mu^{2n+1}}{n!\, (k-1)!\, (2n+1)}\, \left( \frac{\epsilon}{1-\epsilon} \right)^{2n+1-k}\ \Big(\frac 1t\, \int_0^t\, s (Me_\nu)(s)\, (t-s)^{k-1}\, ds \ ,\ e_j\Big)
\end{aligned}
\end{equation}
where in the last summation we have used the $L^2$ convergence of the series. Concerning the last terms we use \cite[eq. 7.415 p. 810]{GR} to write
\[
\frac 1t\, \int_0^t\, s (Me_\nu)(s)\, (t-s)^{k-1}\, ds = t^k\, \int_0^1\, s \, (1-s)^{k-1}\, e^{-st}\, e_\nu(st) ds = \frac{\nu+1}{k(k+1)} \, t^k\, \hyper{1}{1}(\nu+2,k+2,-t)
\] 
where $\hyper{1}{1}$ is the standard confluent hypergeometric function. Moreover, by \eqref{laguerre} and \cite[eq. 7.621(4) p. 822]{GR}
\[
\begin{aligned}
\Big( t^k\, \hyper{1}{1}(\nu+2,k+2,-t)\ ,\ e_j(t) \Big) & = \sum_{m=0}^j\, \comb{j+1}{j-m}\, \frac{(-1)^m}{m!}\, \int_0^\infty\, t^{k+m+1}\, e^{-t}\, \hyper{1}{1}(\nu+2,k+2,-t)\, dt = \\[0.2cm]
& = \sum_{m=0}^j\, \comb{j+1}{j-m}\, \frac{(-1)^m\, (k+m+1)!}{m!}\, \hyper{2}{1}(\nu+2, k+m+2; k+2; -1)
\end{aligned}
\]
where $\hyper{2}{1}$ is the hypergeometric function. Using the previous equations in \eqref{final}, we get the following explicit expression for the general term $a^\mu_{j\nu}$ of the matrix $A_\mu$ defined in \eqref{eig-mat} and \eqref{ar}, which represents the transfer operator $\tilde \PP^{op}_\mu$ on $L^2(m)$
\begin{equation}\label{real-final}
\begin{aligned}
& a^\mu_{j\nu} = \Big(1 - \mathrm{Erf}\Big(\frac{\mu\, \epsilon}{1-\epsilon}\Big)\Big)\, \comb{\nu+j+1}{\nu}\, \frac{1}{2^{\nu+j+3}} + \sum_{\ell=0}^{\nu}(-1)^{\ell}\, \comb{\nu+1}{\nu-\ell}\, \comb{\ell+j+1}{l}\, \frac{1}{2^{\ell+j+2}} +\\[0.2cm]
& + \frac{\nu+1}{(j+1)\sqrt{\pi}}\, \sum_{n=0}^\infty\, \sum_{k=1}^{2n+1}\, \sum_{m=0}^j\,  \comb{2n+1}{k} \comb{j+1}{j-m} \comb{k+m+1}{m} \frac{(-1)^{n+k+m-1}\, \mu^{2n+1}}{n!\, (2n+1)}\, \left( \frac{\epsilon}{1-\epsilon} \right)^{2n+1-k} \cdot \\[0.2cm]
& \cdot \hyper{2}{1}(\nu+2, k+m+2; k+2; -1)
\end{aligned}
\end{equation}
We recall that the parameter $\mu$ determines the level of approximation of the real hole $\tilde H = [0,\frac{\epsilon}{1-\epsilon})$ used in \eqref{to-hole-new}. We need to comment on this approximation in order to compare our numerical results with results by previous papers, and in particular those in \cite{KM}. 

To estimate the scaling of the escape rate $\gamma$ defined in \eqref{escape}, as the hole $H=[0,\epsilon)$, and hence the hole  $\tilde H = [0,\frac{\epsilon}{1-\epsilon})$, shrinks to the empty set, we compute the principal eigenvalue of north-west corner approximations of the matrix $A_\mu$ as $\epsilon$ decreases to 0. We remark that as $\epsilon \to 0^+$, our approximated hole, that is the function $\xi_\mu(x,\frac{\epsilon}{1-\epsilon})$, converges to the function
\[
\xi_\mu(x,0) = \frac 12 - \frac 12 \mathrm{Erf}(\mu\, x)\, ,
\]
which, as $\mu\to +\infty$, converges to zero for all $x>0$. However, if we observe that the appropriate quantity with which the escape rate should be compared is the measure of the hole, we have to look at the integral of $\xi_\mu(x,0)$ on $[0,1]$, and more precisely on $\R^+$ since our approximated holes have effect on the whole real positive axis. Unfortunately the integral shows slow convergence to 0 as $\mu$ diverges. For $\mu=7$, the integral is of order $10^{-2}$, which implies a significant perturbation on the transfer operator. Moreover, a value of $\mu$ greater than 4 in \eqref{real-final} implies that, to have good numerical results, one should consider too many terms in the series on $n$.

\begin{figure}[h]
\begin{center}
\subfigure[]
    {\includegraphics[width=6cm]{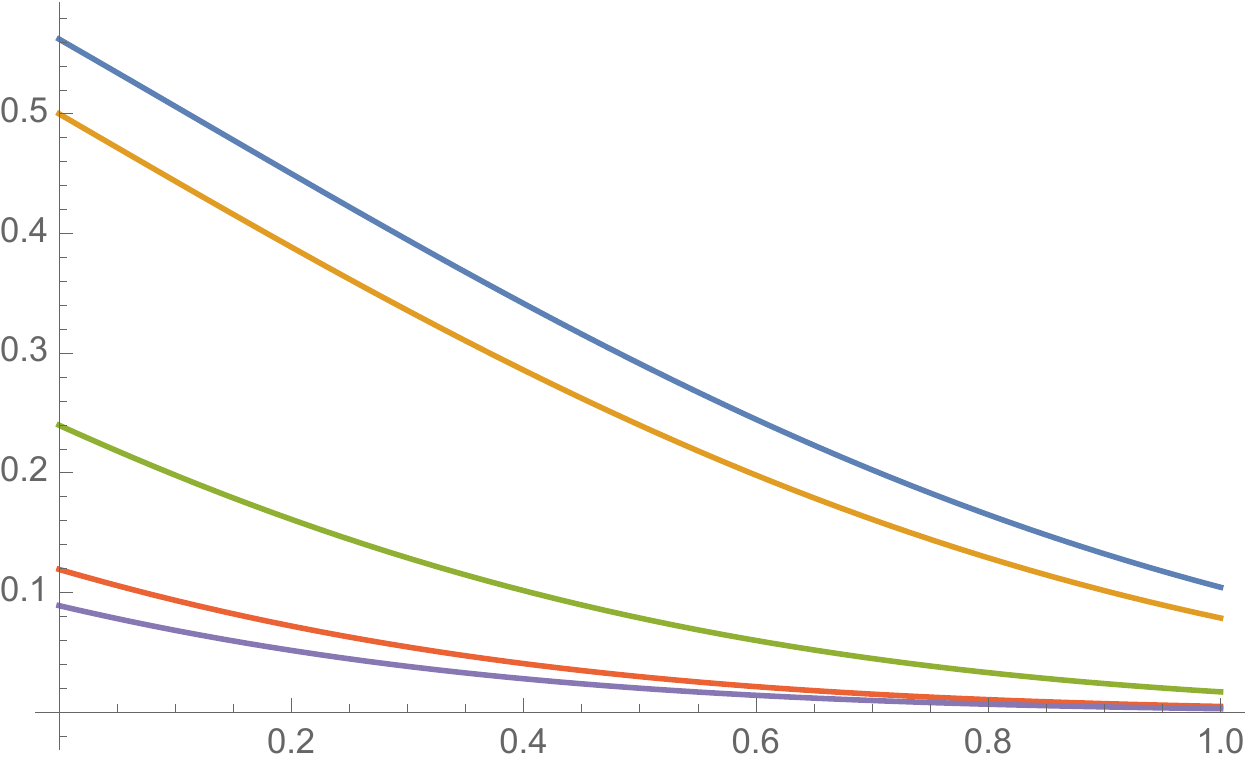}}
    \hspace{0.7cm}
    \subfigure[]
    {\includegraphics[width=6cm]{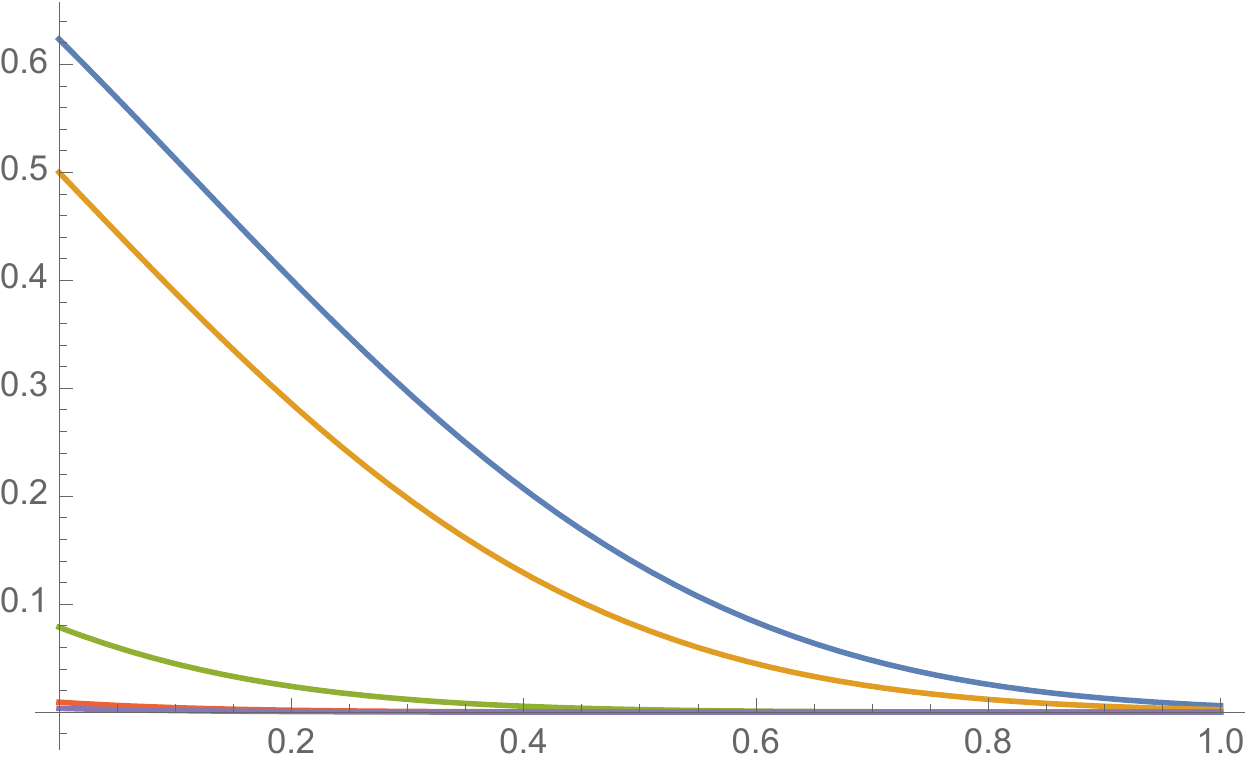}}
\caption{(a) The functions $\xi_\mu(x,\frac{\epsilon}{1-\epsilon})$ for $\mu=1$ and $\epsilon = 0.1, 0, -1, -5, -20$.  (b) The same as (a) for $\mu=2$.} \label{graph-2}
\end{center}
\end{figure}

For this reason, we have changed our point of view, keeping in mind that what is important is to have a sequence of approximated holes with vanishing measure on the real positive axis. This can be obtained for small values of $\mu$, by letting $\epsilon$ diverge to $-\infty$ in $\xi_\mu(x,\frac{\epsilon}{1-\epsilon})$. In Figure \ref{graph-2} we show the behavior of $\xi_\mu(x,\frac{\epsilon}{1-\epsilon})$ on $[0,1]$ for $\epsilon = 0.1, 0, -1, -5, -20$ for $\mu=1$ on the left, and for $\mu=2$ on the right. We notice that the functions decrease on the positive real axis more rapidly for bigger $\mu$, only the first three cases are non-negligible for $\mu=2$, and if we measure their integral on $\R^+$ we have that
\[
\lim_{\epsilon\to -\infty} \int_0^{+\infty}\, \xi_\mu\Big(x,\frac{\epsilon}{1-\epsilon}\Big) \, dx \approx \left\{ \begin{array}{ll} 0.025 & \text{for }\ \mu=1\\ 2\times 10^{-4} & \text{for }\ \mu=2\\  10^{-5} & \text{for }\ \mu=2.5\\  5\times 10^{-7} & \text{for }\ \mu=3\\  10^{-8} & \text{for }\ \mu=3.5 \end{array} \right.
\]
and the integral is already $\approx 10^{-6}$ for $\mu=3$ and $\epsilon = -20$, and $\approx 10^{-7}$ for $\mu=3.5$ and $\epsilon = -10$. 

\begin{figure}[h]
\begin{center}
\includegraphics[width=8cm]{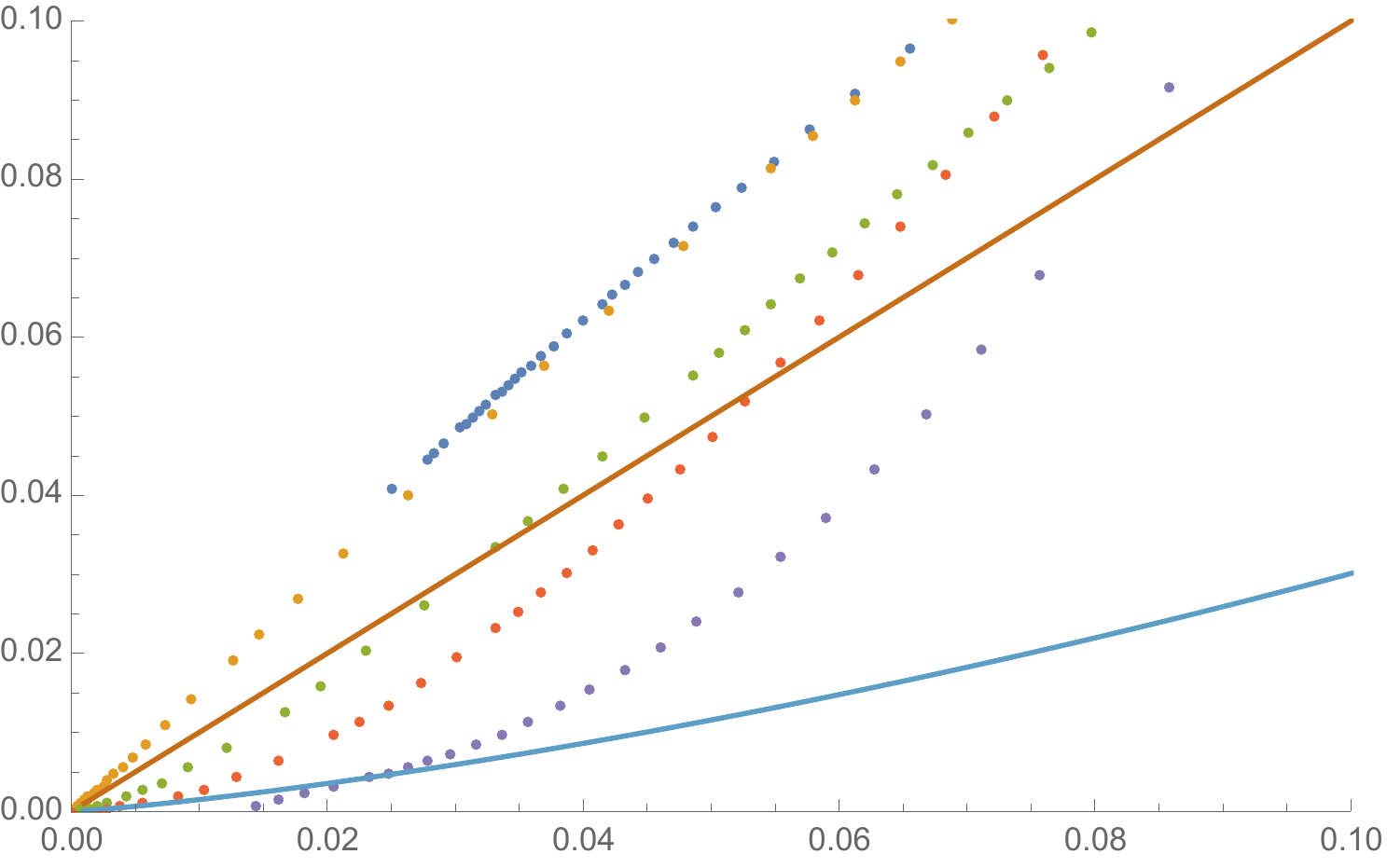}
\caption{The solid lines are the identity and the function $f(t) = \frac{t}{-\log t}$. The dotted lines are the points $(M_\mu(\epsilon),\gamma_\mu(\epsilon))$ for $\mu=1,2,2.5,3$ and 3.5 from the biggest to the lowest.} \label{graph-3}
\end{center}
\end{figure}

Hence we have computed the principal eigenvalue $\lambda_\mu(\epsilon)$ of north-west corner approximations of the matrix $A_\mu$ as $\epsilon$ decreases to $-\infty$, and the principal eigenvalue $\lambda_\infty$ of the same approximations of the matrix associated to the transfer operator of the Farey map without hole. To find the scaling of the escape rate, we have plotted $\gamma_\mu(\epsilon):= -\log (\lambda_\mu(\epsilon) / \lambda_\infty)$ against $M_\mu(\epsilon):= \int_0^{+\infty}\, \xi_\mu\Big(x,\frac{\epsilon}{1-\epsilon}\Big) \, dx$. The results are shown in Figure \ref{graph-3}. The solid lines are the identity and the function $f(t) = \frac{t}{-\log t}$. The dotted lines are the plots of the points $(M_\mu(\epsilon),\gamma_\mu(\epsilon))$ for $\mu=1,2,2.5,3$ and 3.5 from the biggest to the lowest. Notice that for $\mu=1$ the dots stop far from the origin because $M_1(-\infty) \approx 0.025$.

Figure \ref{graph-3} shows that the scaling of the escape rate for the Farey map with shrinking approximated holes is dependent on the shape and the goodness of the approximation. However, as the approximation gets better, that is for $\mu$ big in our case, we find a scaling 
\[
\gamma_\mu(\epsilon) \approx \frac{M_\mu(\epsilon)}{-\log M_\mu(\epsilon)}\qquad \text{as}\quad M_\mu(\epsilon) \to 0^+\, ,
\]
which is consistent with \eqref{kn-mu}, the theoretical result for the Markov approximation of the Farey map studied in \cite{KM} with holes generated by the Markov partition.

\section{Conclusions}

In this paper we have studied the escape rate for the Farey map, an infinite measure preserving system, with a hole including the indifferent fixed point. To our knowledge this is the first time this problem is studied for general holes, since previous results only considered piecewise linear approximations of the map with holes generated by the associated partition.

The problem we consider poses theoretical difficulties in the application of the standard methods for the study of the escape rate, in particular the transfer operator approach, due to the properties of the transfer operator of the Farey map. For this reason, we propose to modify the standard approach to open systems by considering approximations of the hole, by means of a family of functions converging to the indicatrix function of the hole. We believe that this is the main contribution of this paper. Here we have used the family \eqref{famiglia}, and proved in Theorem \ref{funzionale} that the associated transfer operator has the functional analytic properties necessary to pursue the approach already used in \cite{BGI,pap1,pap2} to study the spectral properties of the transfer operator of the Farey map.

Our numerical results suggest that the behavior of the escape rate is indeed dependent on the chosen approximation of the hole, but for functions close to the indicatrix function, we find the same behavior proved in \cite{KM} for the piecewise linear approximation of the map.

\appendix

\section{Proof of Theorem \ref{funzionale}}
From \eqref{emme-tilde}, we have to compute the inverse Laplace transform of the function $\xi$. For this we use the series expansion in \eqref{series} and write
\[
\begin{aligned}
\LL^{-1}\left[\xi_\mu\left(\frac 1x, a\right)\right](t) & = \LL^{-1}\left[ \frac 12 - \frac 12\text{Erf}\left(\mu \, \left(\frac 1x  - a\right)\right)\right] = \\[0.2cm]
= & \frac 12 \, \delta_0(t) - \frac 12\, \LL^{-1}\left[ \frac{2}{\sqrt{\pi}}\, \sum_{n=0}^\infty\, \frac{(-1)^n\, \mu^{2n+1}}{n!\, (2n+1)} \left(\frac 1x  - a\right)^{2n+1}\right](t) = \\[0.2cm]
= & \frac 12 \, \delta_0(t) - \frac{1}{\sqrt{\pi}}\, \sum_{n=0}^\infty\, \frac{(-1)^n\, \mu^{2n+1}}{n!\, (2n+1)}\, \LL^{-1}\left[ \sum_{k=0}^{2n+1}\, \comb{2n+1}{k} \frac{1}{x^k}{(-a)^{2n+1-k}}\right](t) = \\[0.2cm]
= & \frac 12 \, \delta_0(t) - \frac{1}{\sqrt{\pi}}\, \sum_{n=0}^\infty\, \frac{(-1)^{n-1}\, (\mu a)^{2n+1}}{n!\, (2n+1)}\, \delta_0(t) + \\[0.2cm]
& - \frac{1}{\sqrt{\pi}}\, \sum_{n=0}^\infty\, \sum_{k=1}^{2n+1}\, \comb{2n+1}{k} \, \frac{(-1)^{n+k-1}\, \mu^{2n+1}\, a^{2n+1-k}}{n!\, (2n+1)}\, \frac{t^{k-1}}{(k-1)!} = \\[0.2cm]
= & \frac 12 \Big( 1+ \text{Erf}(\mu a) \Big)\, \delta_0(t) - \frac{1}{\sqrt{\pi}}\, \sum_{n=0}^\infty\, \sum_{k=1}^{2n+1}\, \comb{2n+1}{k} \, \frac{(-1)^{n+k-1}\, \mu^{2n+1}\, a^{2n+1-k}}{n!\, (2n+1)}\, \frac{t^{k-1}}{(k-1)!}\, ,
\end{aligned}
\]
where $\delta_0$ denotes the Dirac delta at the origin, and we have used the standard results
\[
\LL^{-1}[x^{-k}](t) = \frac{t^{k-1}}{(k-1)!} \quad \text{for}\ t\in \R^+\, , \, k>0\, , \qquad \LL^{-1}[1](t) = \delta_0(t)\, .
\]
Moreover, from the second to the third line, we have used that
\[
\int_0^\infty\, e^{-xt}\, \left[\sum_{n=0}^\infty\, \frac{(-1)^n\, \mu^{2n+1}}{n!\, (2n+1)} \left( \delta_0(t) + \sum_{k=1}^{2n+1}\, \comb{2n+1}{k} \, \frac{(-a)^{2n+1-k}}{(k-1)!}\, t^{k-1} \right)\right] \, dt =
\]
\[
= \sum_{n=0}^\infty\, \frac{(-1)^n\, \mu^{2n+1}}{n!\, (2n+1)}\, \int_0^\infty\, e^{-xt}\, \left( \delta_0(t) + \sum_{k=1}^{2n+1}\, \comb{2n+1}{k} \, \frac{(-a)^{2n+1-k}}{(k-1)!}\, t^{k-1} \right) \, dt 
\]
for $x$ big enough. Hence, by the analytical continuation of the Laplace transform, the equality is verified for all $x>0$. It follows that the passage of the $\LL^{-1}$ operator inside the summation, from the second to the third line, is justified.

Since 
\[
\int_0^t\, s(M\varphi)(s) \delta_0(t-s)\, ds = t (M\varphi)(t)
\]
from \eqref{emme-tilde} we obtain the first term on the right hand side of \eqref{trad-op}. To finish we need to show that
\begin{equation}\label{stima-l2}
\begin{aligned}
& \frac 1t\, \int_0^t\, s\, (M\varphi) (s) \, \sum_{n=0}^\infty\, \sum_{k=1}^{2n+1}\, \comb{2n+1}{k} \, \frac{(-1)^{n+k-1}\, \mu^{2n+1}\, a^{2n+1-k}}{n!\, (2n+1)}\, \frac{(t-s)^{k-1}}{(k-1)!}\, ds =\\[0.2cm]
& = \sum_{n=0}^\infty\, \sum_{k=1}^{2n+1}\, \comb{2n+1}{k} \, \frac{(-1)^{n+k-1}\, \mu^{2n+1}\, a^{2n+1-k}}{n!\, (2n+1)}\, \frac 1t\, \int_0^t\, s\, (M\varphi) (s) \, \frac{(t-s)^{k-1}}{(k-1)!}\, ds
\end{aligned}
\end{equation}
First of all we recall from \cite[eq. 13.2.2 p. 322 and eq. 13.6.19 p. 328]{nist} that for any $\varphi \in L^2(m)$ we have
\[
\sum_{s=0}^{2n}\, \frac{(-2n)_s}{(2)_s\, s!}\, \Big(\frac ta\Big)^s = \hyper{1}{1}\Big(-2n,2;\frac ta\Big) = \frac{1}{2n+1}\, e_{2n} \Big(\frac ta\Big)
\]
where $(k)_s = k (k+1)\dots(k+s-1)$ is the Pochhammer symbol, $\hyper{1}{1}$ is the standard confluent hypergeometric function, and $e_\nu$ is the Laguerre polynomial defined in \eqref{laguerre}. Hence
\[
\sum_{n=0}^\infty\, \sum_{k=1}^{2n+1}\, \comb{2n+1}{k} \, \frac{(-1)^{n+k-1}\, \mu^{2n+1}\, a^{2n+1-k}}{n!\, (2n+1)}\, \frac{(t-s)^{k-1}}{(k-1)!} = \sum_{n=0}^\infty\, \frac{(-1)^{n}\, \mu^{2n+1}\, a^{2n}}{n!\, (2n+1)}\, e_{2n} \Big(\frac{t-s}{a}\Big)
\]
We are then reduced to study the increase of the terms
\[
\Big\| \frac 1t\, \int_0^t\, s\, (M\varphi) (s)\, e_{2n} \Big(\frac{t-s}{a}\Big) \, ds \Big\|_{L^2(m)}^2 = \int_0^\infty\, \frac 1t\, e^{-t}\, \left(\int_0^t\, s\, e^{-s} \, \varphi (s)\, e_{2n} \Big(\frac{t-s}{a}\Big) \, ds\right)^2\, dt
\]
Standard manipulations show that
\[
\begin{aligned}
& \Big\| \frac 1t\, \int_0^t\, s\, (M\varphi) (s)\, e_{2n} \Big(\frac{t-s}{a}\Big) \, ds \Big\|_{L^2(m)}^2 \le \| \varphi\|_{L^2(m)}^2\, \int_0^\infty\, \frac 1t\, \, e^{-t}\, \int_0^t\, s\, e^{-s} \, e_{2n}^2 \Big(\frac{t-s}{a}\Big) \, ds\ dt  = \\[0.2cm]
& = \| \varphi\|_{L^2(m)}^2\, \int_0^\infty\, e^{-u}\, \int_0^u\, \Big(\frac{u-v}{u+v}\Big)\, e_{2n}^2 \Big(\frac{v}{a}\Big) \, dv\ du \le a\, \| \varphi\|_{L^2(m)}^2\, \int_0^\infty\, e^{-u}\, \int_0^{u/a}\, e_{2n}^2(v)\, dv\, du\, .
\end{aligned}
\]
By \eqref{laguerre}, we have
\[
e_{2n}^2 (v) = \sum_{i,j=0}^{2n}\, \comb{2n+1}{2n-i} \comb{2n+1}{2n-j} \frac{(-1)^{i+j}}{i!\, j!}\, v^{i+j}
\]
hence
\[
\begin{aligned}
& \int_0^\infty\, e^{-u}\, \int_0^{u/a}\, e_{2n}^2(v)\, dv\, du = \sum_{i,j=0}^{2n}\, \comb{2n+1}{2n-i} \comb{2n+1}{2n-j} \frac{(-1)^{i+j}}{i!\, j!}\, \int_0^\infty\, e^{-u}\, \int_0^{u/a}\, v^{i+j}\, dv\, du = \\[0.2cm]
& = \sum_{i,j=0}^{2n}\, \comb{2n+1}{2n-i} \comb{2n+1}{2n-j} \frac{(-1)^{i+j}}{a^{i+j+1}}\, \frac{1}{i!\, j!\, (i+j+1)}  \int_0^\infty\, e^{-u}\, u^{i+j+1}\, du = \\[0.2cm]
& = \sum_{i,j=0}^{2n}\, \comb{2n+1}{2n-i} \comb{2n+1}{2n-j} \frac{(-1)^{i+j}}{a^{i+j+1}}\, \comb{i+j}{i}
\end{aligned}
\]
Using the very crude estimate $\comb{k}{h}\le 2^k$ for all $h=0,\dots,k$, and $i+j\le 4n$, we obtain
\[
\int_0^\infty\, e^{-u}\, \int_0^{u/a}\, e_{2n}^2(v)\, dv\, du \le (2n+1)^2\, 2^{4(2n+1)}\, \max\{ a^{-s}\, :\, s=1,\dots, 4n+1\}\, .
\]
Hence
\[
\sum_{n=0}^\infty\, \frac{(-1)^{n}\, \mu^{2n+1}\, a^{2n}}{n!\, (2n+1)}\, \Big\| \frac 1t\, \int_0^t\, s\, (M\varphi) (s)\, e_{2n} \Big(\frac{t-s}{a}\Big) \, ds \Big\|_{L^2(m)} \le \| \varphi\|_{L^2(m)}\, \sum_{n=0}^\infty\, \frac{(2n+1)\, c^{2n+1}}{n!}
\]
where
\[
c = \left\{ \begin{array}{ll} 16\, \frac{\mu}{a}\, , & \text{if }\, |a|\le 1\\[0.2cm] 16\, \mu\, a\, , & \text{if }\, |a|> 1 \end{array} \right.
\]
In any case we obtain total convergence in $L^2(m)$ for the right hand side of \eqref{stima-l2}, hence \eqref{stima-l2} holds in the $L^2$-sense.

This proves \eqref{trad-op} and, together with the boundedness of $M$ on $L^2(m)$, shows that $\tilde M$ is a bounded operator. This finishes the proof of Theorem \ref{funzionale}.

\end{document}